\documentclass[11pt]{amsart}

\scrollmode
\usepackage{amssymb,amsmath}
\usepackage{graphicx}
\usepackage{datetime}

\theoremstyle{plain}
\newtheorem{thm}{Theorem}[section]

\newtheorem{rem}[thm]{Remark}

\newtheorem{ass}[thm]{Assumption}

\numberwithin{equation}{section}

\def\svdots{\vbox{\baselineskip=1.5pt\lineskiplimit=0pt
	\kern1.5pt \hbox{$\ss .$}\hbox{$\ss .$}\hbox{$\ss .$}}}

\def\sddots{\mathinner{\raise3pt\vbox{\hbox{$\ss .$}}
    \raise1.5pt\hbox{$\ss .$}\hbox{$\ss .$}}}

\newcommand{\R}{{\mathbb R}}

\newcommand{\Rn}{{\R^n}}

\newcommand{\Rnxn}{{\R^{n\times n}}}

\newcommand{\cont}{{\mathcal C}}

\newcommand{\itext}[1]{{\qquad\text{#1}\qquad}}
\newcommand{\abs}[1]{\left| #1 \right|}

\newcommand{\smat}[1]{ \left[\begin{smallmatrix} #1 \end{smallmatrix}\right]}

\newcommand{\bmat}[1]{ \begin{bmatrix}#1\end{bmatrix}}

\let\hat\widehat

\renewcommand{\ss}{\scriptstyle}

\def\BSigma{{\boldsymbol{\Sigma}}}

\renewcommand{\vec}[1]{{\boldsymbol{#1}}}

\def\s{\vec{s}}

\def\x{\vec{x}}
\def\y{\vec{y}}

\def\yes{Yes}
\newcommand{\ForUsOnly}[1]{\ifx\UsOnly\yes {\em #1}\fi}

\newtheorem{theorem}{Theorem}

\newtheorem{lemma}[theorem]{Lemma}

\newtheorem{exam}[theorem]{Example}

\theoremstyle{definition}
\newtheorem{definition}{Definition}

\begin{document}

\title{Master stability function for  \\ piecewise smooth networks}

\author{Luca Dieci}
\address{School of Mathematics, Georgia Tech, Atlanta, GA 30332, USA}
\email{dieci@math.gatech.edu}
\thanks{}
%    author one information
\author{Cinzia Elia}
\address{Dipartimento di Matematica, Univ. of Bari, I-70100, Bari, Italy}
\email{cinzia.elia@uniba.it}
\thanks{}
%    author two information

%\date{\today}

\dedicatory{}

\subjclass[2010]{Primary 34A36, Secondary 34D06, 34D08}

\keywords{Piecewise smooth networks, synchronization, fundamental matrix solution,
master stability function, Floquet mutlipliers, Floquet exponents}

\null\hfill Version of \today $, \,\,\,$ \xxivtime

\begin{abstract}
We consider a network of identical piecewise smooth systems that synchronizes on the
manifold given by a periodic orbit of a single agent.
We explicitly characterize the fundamental matrix solution of the network along the synchronous 
solution and extend the {\sl Master Stability Function} tool to the present case of 
non-smooth dynamics.
\end{abstract}

\maketitle

\pagestyle{myheadings}
\thispagestyle{plain}

\pagestyle{myheadings}
\thispagestyle{plain}
\markboth{L Dieci, C Elia}{MSF for PWS networks}

\section{Introduction}
Synchronization of dynamical networks is a fascinating, widely studied, and impactful phenomenon; e.g., 
see \cite{Winfree, Buck, Kuramoto} for early applications in the
applied sciences, and the review
\cite{Dorfler-Bullo} --and the many references there-- for a thorough
account on the topic of oscillators synchronization.
In the most typical and studied circumstance, one has a network of $N$ nodes,
the so-called {\it agents}, that obey $N$ identical nonlinear differential equations with vector field $f$, 
coupled through linear anti-symmetric coupling.  
The key concern is to find conditions that tie together the strength of the coupling, 
the structure of the network, and the agent's dynamics, in such a way that the network synchronizes.
This problem has been extensively studied, under a number of different scenarios, for
networks of smooth systems; e.g., see the work of Pecora and coauthors \cite{BarahonaPecora, Pecora.Carroll} for a
study exploiting Lyapunov exponents ideas, and see \cite{Hale} for a study more 
along the lines of the theory of dissipative attractors. 

After the cited works of Pecora and coauthors, probably the most widely adopted
and successful tools to infer convergence to, and/or stability of,
a synchronized solution in networks of smooth dynamical systems has been that of the
master stability function (MSF).  However,
as remarked in \cite{CdLdB}, this {\sl ``approach requires some degree
of smoothness in the agents’ vector fields ... and extensions
need to be found''} when dealing with piecewise smooth systems.
Our goal in this work is to provide such extension.

In fact, in the present work, we study the synchronization problem for networks where each agent satisfies a nonlinear piecewise smooth
system, that is the underlying dynamical system of each agent is governed by a vector field $f$ which is
piecewise smooth.  We will call these  {\sl piecewise smooth networks}.  
That is, we consider the network 
\begin{equation}\label{PWSnetwork1}
	\dot {x_i}=f({x_i})+\sigma \sum_{j=1}^N a_{ij}E(x_j-x_i),\quad x_i\in \R^n,\,\,\
	i=1\dots, N\ ,
\end{equation}  
where $f$ is only piecewise smooth:
$f(x_i)=\left \{ \begin{matrix} 
	f^+(x_i), & h(x_i)>0 \\ f^-(x_i) & h(x_i)<0 \end{matrix} \right .$, 
$i=1, \ldots, N$, and $h(x)=0$ describes the discontinuity manifold.
In \eqref{PWSnetwork1},
$A\in \R^{N\times N}$ is the adjacency matrix of the 
graph describing the network (assumed to be 
undirected, simple and connected, so that $A$ is symmetric), $\sigma \ge 0$ is the the coupling strength, 
and $E\in \Rnxn$ is the coupling matrix describing which components
of the two agents $x_i$ and $x_j$, $i\ne j$, are connected to one another.
Next, we let $D$ be the diagonal matrix with elements $d_{ii}=\sum_{j=1}^N a_{ij}$, and  
let $L=-D+A$, $L \in \R^{N \times N}$, the graph Laplacian. Then,
using Kronecker product notation, we rewrite
 \eqref{PWSnetwork1} as
\begin{equation}\label{PWSnetwork}\begin{split}
\dot {\x}&=F({\x})+\sigma M{\x},\itext{where} \\
{\x}&=\bmat{x_1 \\ \vdots \\ x_N}\in \R^{nN}, \,
F({\x})=\bmat{f(x_1)\\ \vdots \\ f(x_N)} \in \R^{nN}, 
\, M= L \otimes E \in \R^{nN\times nN}\ .
\end{split}\end{equation}  
As noted, our present interest is in the case when, taking $\sigma=0$ in \eqref{PWSnetwork}, 
each agent obeys identical piecewise smooth (PWS) dynamics:
\begin{equation}\label{nonsmooth_eq}
	\dot x_i=f(x_i)=\left \{ \begin{matrix} 
		f^+(x_i), & h(x _i)>0 \\ f^-(x_i), & h(x_i)<0 \end{matrix} \right ., 
	\qquad i=1, \ldots, N,
\end{equation}
with $f^\pm: \R^n \to \R^n $ smooth vector fields, and $h: \Rn\to \R$ is
assumed to be at least $\cont^2$.
For each agent, the manifold of discontinuity is the zero set $\{x\in \Rn \ :\,\ h(x)=0\}$,
and we will use the following notation:
\begin{equation}\label{notations_eq}
	\Sigma=\{x \in \R^n \,\ | \,\ h(x) =0\},\,\
	R^\pm=\{x \in \R^n \,\ | \,\ h(x) \gtrless 0\}.
\end{equation}
\begin{rem}
A typical case we have seen in many applications is to have	
$h(x)=c^Tx-b$, but this is not necessary in our work.
\end{rem}
As customary, we say that a point $x \in \Sigma$ is a \emph{transversal crossing point} if 
\begin{equation}\label{CrossPt}
(\nabla h(x)^Tf^-(x) )(\nabla h(x)^T f^+(x) )>0,
\end{equation}
and it is an \emph{attractive sliding point} if 
\begin{equation}\label{SlidePt}
\nabla h(x)^Tf^-(x)>0, \quad \nabla h(x)^T f^+(x)<0.
\end{equation}
On $\Sigma$, sliding will be assumed to take place in the sense of Filippov, 
whereby on $\Sigma$ the dynamics of an agent is given by
\begin{equation}\label{SlidingDyn}
	\dot x =f_\Sigma:=(1-\alpha)f^{-}(x)+\alpha f^{+}(x)\ ,\quad
	\alpha=\frac{\nabla h(x)^Tf^-(x)}{\nabla h(x)^T(f^-(x)-f^+(x))}\ .%	\quad i=1, \ldots, N.
\end{equation}
Finally, a point 
$\bar x \in \Sigma$ is called \emph{tangential exit point into} $R^-$ if a trajectory $x(t)$
sliding on $\Sigma$ reaches it at some value $\bar t$ and there it holds that
\[\nabla h(\bar x)^Tf^-(\bar x)=0, \quad \nabla h(\bar x)^Tf^+(\bar x)<0, \quad
\left[\frac{d}{dt}\nabla h(x(t))^Tf^-(x(t))\right]_{t=\bar t}<0,  \]
and similarly for a tangential exit point into $R^+$.  The combination of transversal crossings,
transversal entries on $\Sigma$, and tangential exits from $\Sigma$, are called {\sl generic events},
or simply {\sl events}.

Let the single agent \eqref{nonsmooth_eq} have a limit cycle
with a finite number of events,
and not entirely contained in $\Sigma$ (see
\cite{Galvanetto_95}, \cite{Dieci.Elia_2014}, \cite{Coombes.Thul_2016} for examples of self sustained 
oscillations in discontinuous systems with partial sliding along the discontinuity manifold).  Let
$x_S(t)$ be the corresponding $T$-periodic solution.  Then, given the structure of $M$, 
the function ${\x_S}(t)=\smat{x_S(t)\\ \svdots\\  x_S(t)}$ is a periodic solution of 
\eqref{PWSnetwork} of period $T$; we will call this the
{\sl synchronized manifold} or simply the {\sl synchronous solution}.
However, even if ${ x_S}(t)$ happened to be asymptotically stable for the single agent, there is no 
guarantee that $\x_S$ be stable for the network dynamics for all values of $\sigma$; 
further, when $N$ is large, the numerical study
of the stability of $\x_S(t)$ may be prohibitively expensive. This issue can be overcome by
extending the {\sl Master Stability Function} (MSF) tool of Pecora and Carroll,
see \cite{Pecora.Carroll}, to PWS networks.

The MSF technique relies on exploiting the structure of the fundamental matrix solution of the 
network, and for this reason in the present work our goal is two-fold.
When
the network synchronizes on ${\x_S}(t)$,
first we will give the explicit expression of the fundamental matrix
solution along the synchronized manifold. Then, we will extend the {\sl Master Stability Function} (MSF) 
to piecewise smooth networks.
Many authors have considered piecewise smooth networks, and some important
studies have been made to resolve the outstanding concern of how
to infer asymptotic convergence in networks of piecewise-smooth systems.
Notable examples are the recent work of \cite{CdLdB2} where discontinuous
diffusive coupling is adopted, and the works  \cite{Coombes.Thul_2016, Coombes.Thul_2018}, where
the authors use the MSF approach 
to study limit cycles in piecewise-linear systems.  
However, a rigorous justification of the {\bf use of the MSF} for general, {\bf nonlinear, PWS networks}
appears to be lacking, and it is our purpose to give it in this work.

A plan of the paper is as follows.  In Section \ref{FMS}, we derive the precise form of
the monodromy matrix along the synchronous solution.  In Section \ref{MSF}, we 
extend the MSF tool to piecewise smooth networks.  Finally, in Section \ref{GalvaSection}
we give detailed numerical study of a network arising in mechanical vibrations and infer that, for
a range of values of $\sigma$, the synchronous manifold is stable. \\
\noindent{\bf Notation}.\\
$e \in \R^N$ is the vector with all elements equal to $1$, so that 
$\x_S=e \otimes x_S$ is the synchronous solution in $\R^{nN}$.\\ 
$h_i({\x})  =  h(x_i)$, $i=1, \ldots, N$.
$\Sigma_i =   \{ {\x} \in \R^{nN} \,\  | \,\ h_i({\x})=0 \}$,  $i=1, \ldots ,N$, and 
${\BSigma} =  \cap_{i=1}^N\Sigma_i$.

\section{Fundamental matrix solution for synchronous periodic solutions}\label{FMS}

The main difficulties we need to address in this section are the following.
\begin{itemize}
\item[(i)]  The network \eqref{PWSnetwork}  has $N$ discontinuity manifolds and solutions might slide on
the intersection of two or more manifolds (in fact, as we will see, a synchronous periodic solution $\x_S$
with $x_S$ having a sliding portion, will necessarily slide on the intersection of all $N$ manifolds).
But, in general, the sliding vector field on the intersection of the discontinuity manifolds is not uniquely 
defined and we need to address how this impacts the form of the fundamental matrix of the linearized system.  
In the specific case we consider here, there is no such ambiguity, see 
Lemma \ref{OnSigma_prop} and Theorem \ref{X_thm}.
\item[ii)] The monodromy matrix along a periodic solution of the piecewise system 
\eqref{nonsmooth_eq} is not continuous: 
it has jumps at the entry points (crossing or sliding) on the discontinuity manifold.  
These jumps are taken into account via so called jump or saltation matrices, whose
scope is to transform the vector field at the entry time, say $t^-$,
into the vector field at the exiting time, $t^+$.
The correct expression for such matrices is well known in the literature in the case of a single discontinuity
manifold (see \cite{Aizerman.Gantmacher}, \cite{Mueller}, \cite{Leine}). 
However, synchronous sliding solutions have to slide on the intersection of $N$ discontinuity 
manifolds and the fundamental matrix solution along a synchronous solution must take into 
account jumps at this intersection. In the literature, there are results about these jump
matrices relative to the intersection of two discontinuity manifolds, see \cite{Ivanov} for the
case of crossing and \cite{Dieci.Lopez_2011} for the case of sliding, but no result exist for
the intersection of more than two manifolds. Surely this must be because, in the case of sliding solutions,
there is no uniquely defined Filippov sliding vector field on the intersection of discontinuity manifolds,
as noted in i) above. However, this is not the only issue.  Indeed, in general, 
on the intersection  of discontinuity manifolds,
the jump matrix itself is not uniquely defined, even if we are willing to select
a specific sliding vector field (again, see \cite{Ivanov} for the case of
crossing and \cite{Dieci.Lopez_2011} for the case of sliding).  This being the case,
the fundamental matrix solution cannot be defined in a unique way. 
Theorems \ref{S14_thm} and \ref{S4Sigma_thm} deal with this aspect in case of the
synchronous periodic solution $\x_S$ of \eqref{PWSnetwork}.
\end{itemize}

After the expression for the monodromy matrix is arrived at, in Section \ref{MSF} we will
see how to extend the MSF tool to PWS networks.
  
For the above reasons, hereafter we derive the monodromy matrix along the synchronous solution 
$\x_S$ of \eqref{PWSnetwork}. 
The main results are given in Theorem \ref{S14_thm} and \ref{S4Sigma_thm}, where we show that the 
saltation matrices can be represented in a unique way. 
Recalling that $\x_S=\smat{x_S\\ \svdots \\ x_S}$, where $x_S$ is the periodic solution of
a single agent \eqref{nonsmooth_eq}, we will assume that $x_S$ has a finite number of generic events.  
Because of this, we will make the following convenient assumption on the dynamics of $x_S$.

\begin{figure}\label{po_types}
\begin{minipage}[c]{0.48\linewidth}
	\centering
	\includegraphics[width=8cm,height=8cm]{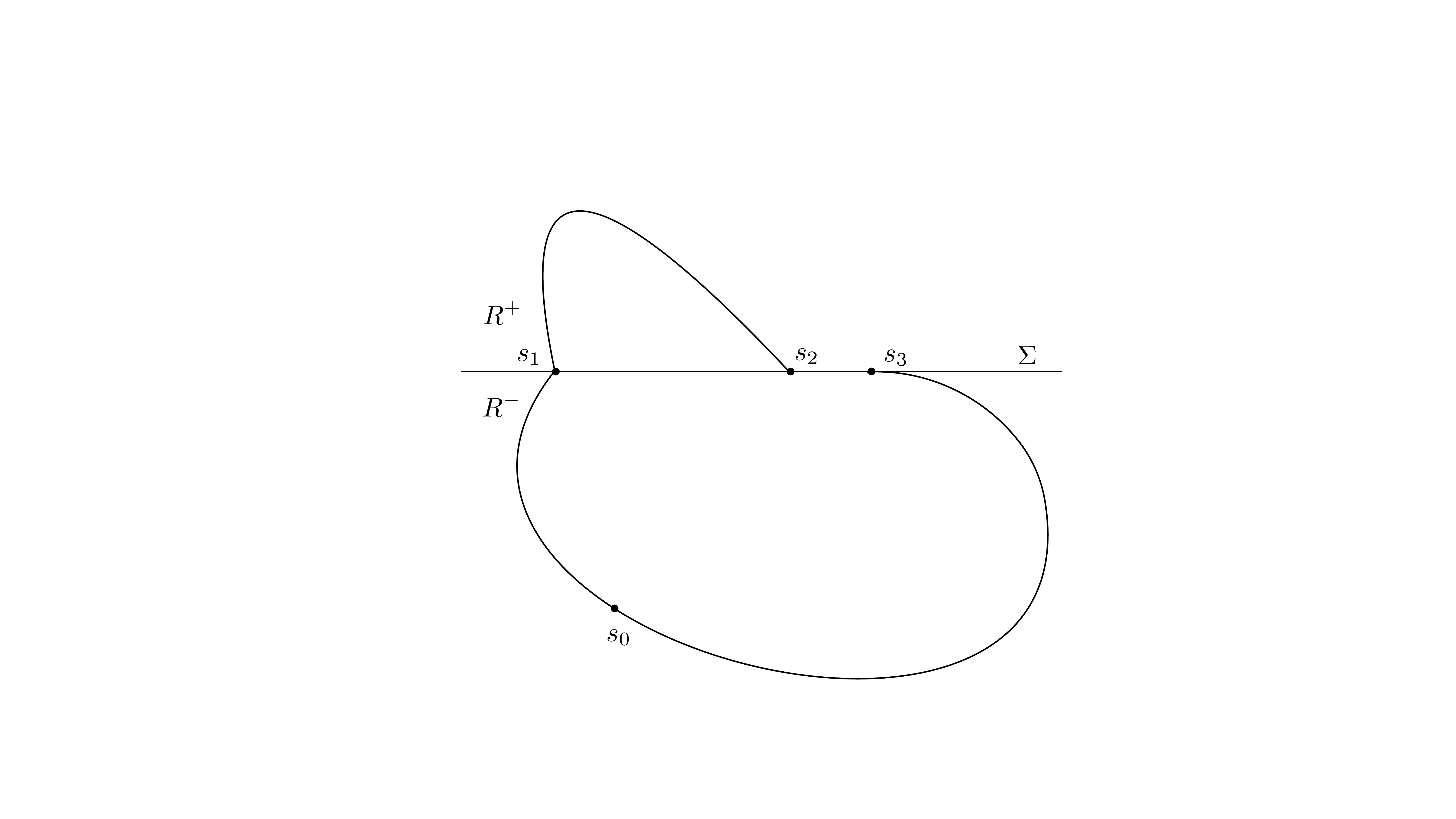}
\end{minipage}
\begin{minipage}[c]{0.48\linewidth}
	\centering
	\includegraphics[width=8cm,height=8cm]{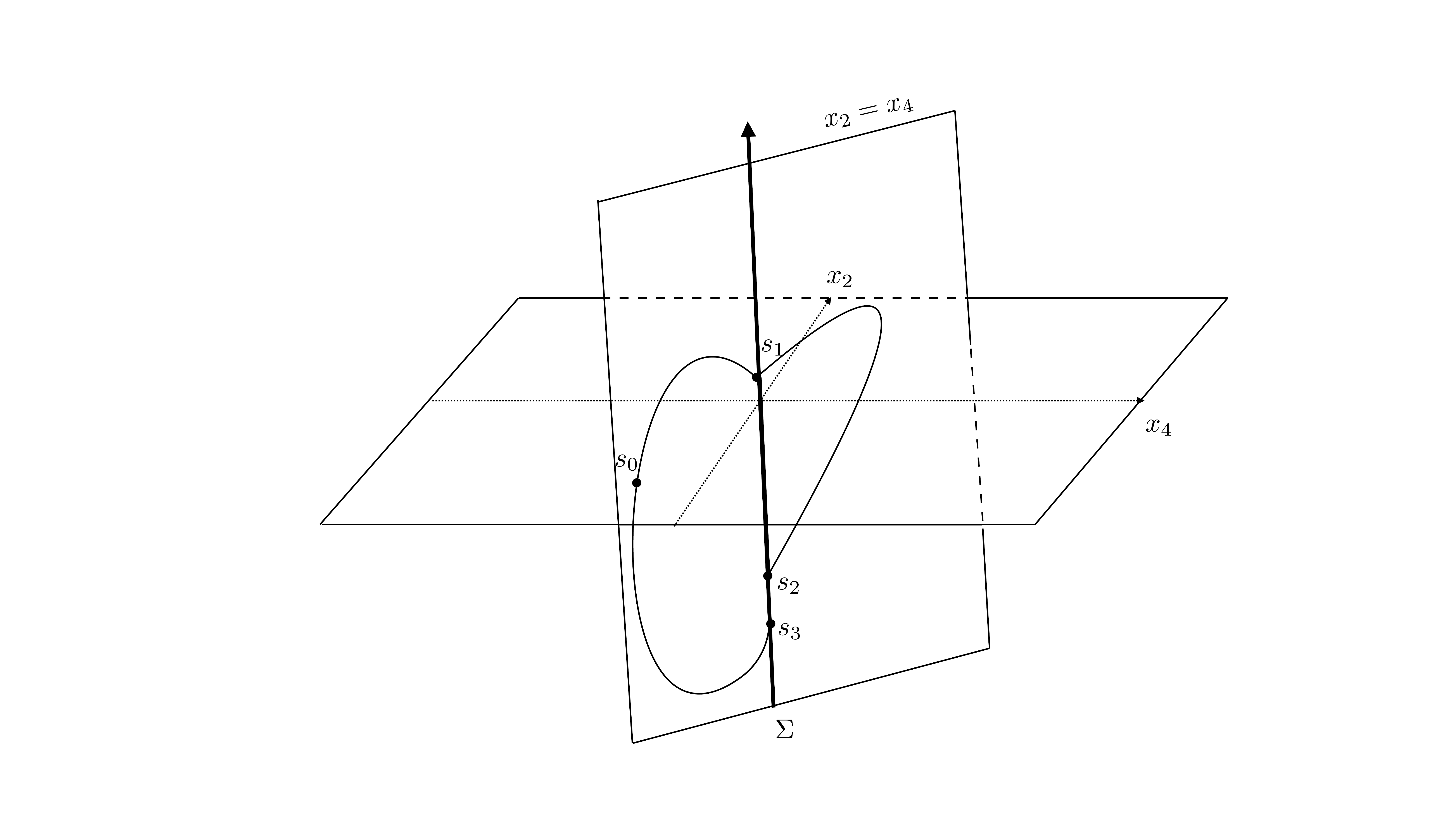}
\end{minipage}
	\caption{Schematic of Assumption \ref{xs_ass}, on the left, and of the periodic orbit
	for	the network, on the right.}
\end{figure}

\begin{ass}\label{xs_ass}
We assume that \eqref{nonsmooth_eq} has a periodic solution $x_S(t)$ that: 
\begin{itemize}
\item[0)] At $t=0$, $x_S(t)=s_0$ is in $R^-$;
\item[1)] At $t=t_1$, $x_S$ crosses $\Sigma$ transversally 
at the point $s_1=x_S(t_1)$ to enter $R^+$; 
\item[2)] At $t=t_2$, $x_S$ reaches transversally the attractive sliding point $s_2=x_S(t_2)\in \Sigma$ 
and $x_S$  begins sliding on $\Sigma$;
\item[3)] At $t=t_3$, $x_S$ reaches the tangential exit point $s_3=x_S(t_3)$, and it
leaves $\Sigma$ to enter into $R^-$;
\item[4)] At $t=T$, $x_S$ is back at $s_0$: $x_S(T)=s_0$. 
\end{itemize}
\end{ass}
A Figure of the above situation is on the left of Figure \ref{po_types}.
\begin{rem}\label{wlog_ass}
Other than the need for a finite number of generic events, the results in this section do not depend on the
particular structure of $x_S(t)$ given in Assumption \ref{xs_ass} and can be immediately 
extended to any finite number of generic crossings, sliding segments, and tangential exits, 
of the periodic orbit of \eqref{nonsmooth_eq}.
\end{rem}

Now, for $N$ agents, 
there are $2^N$ subregions (and corresponding vector fields), and we 
can represent them using a tree diagram with $2^N$ branches. 
We number the regions, and the vector fields, from $1$ to $2^N$ following the branches of the tree.
\begin{exam}\label{regions_ex}
For $N=3$, we have the following correspondence between region numbering and signs
of $h_1$, $h_2$ and $h_3$:
\begin{equation*}
\begin{matrix} 
1& 2& 3& 4& 5& 6 & 7 & 8 \\
(---) & (--+) & (-+-) & (-++) & (+--) & (+-+) & (++-) & (+++) 
\end{matrix}
\end{equation*}  
In each subregion $R_j$, the vector field in \eqref{PWSnetwork} is 
$F(\x)=F_j({\x})=\bmat{f^\pm(x_1)\\ f^\pm(x_2)\\ f^\pm(x_3)}$, $j=1, \ldots, 8$, where 
in $f^\pm$ we select the sign in agreement with the region numbering above. 
For example $F_3({\x})=\bmat{f^-(x_1)\\ f^+(x_2)\\ f^-(x_3)}$.
\end{exam}

\begin{rem}\label{attractive_rem}
It is simple, but important, to observe that 
if $x \in \Sigma$ is an attractive sliding point for the single agent, then 
${\x}=e \otimes x \in {\BSigma}$ is an attractive sliding point on ${\BSigma}$ for the full network.
Indeed, it is immediate to verify that $F_j({\x})$ points toward $\Sigma_i$ for all 
$i=1,\ldots ,N$, and $j=1, \ldots, 2^N$, i.e.,  
\[ \nabla h_i({\x})^TF_j({\x}) > 0, \,\quad
\nabla h_i({\x})^TF_j({\x}) < 0 .
 \]
Similarly, if $x \in \Sigma$ is a tangential exit point into $R^-$ (respectively, $R^+$) for the single agent, 
then ${\x}=e \otimes x \in {\BSigma}$ is a tangential exit point into $R_1$ (respectively, $R_{2^N}$) 
for the full network.
\end{rem}
 
Remark \ref{attractive_rem} justifies the following fact. 
Let $x_S(t)$ satisfy Assumption \ref{xs_ass} and 
let ${\s}_j=e \otimes s_j$, where $s_j$ is defined in Assumption \ref{xs_ass}, $j=1,\ldots, 4$.
Then, the synchronous solution $\x_S(t)=(e \otimes x_S(t))$ obeys the following evolution:
\begin{itemize}
\item[0)] At $t=0$, ${\x_S}(0)={\s}_0$ is in $R_1$;
\item[1)] At $t=t_1$, $\x_S$ crosses $R_1$ at ${\s}_1$ and enters into $R_{2^N}$; 
\item[2)] At $t=t_2$, $\x_S$  reaches the attractive sliding point ${\s}_2$ and starts sliding\footnote{Although, in general, sliding along ${\BSigma}$ is not unambiguously
	defined, presently this is not a concern, since we are just describing the evolution of
	the specific $\x_S$.} 
along ${\BSigma}$;
\item[3)] At $t=t_3$, $\x_S$ exits ${\BSigma}$ at the tangential exit point $x={\s}_3$ and enters 
into $R_1$;
\item[4)] At $t=T$, $\x_S$ reaches ${\s}_0$. 
\end{itemize}
\begin{rem}\label{synch_rem}
Note that the synchronous solution $\x_S$ satisfies the following : 
i) it can only evolve in the regions $R_1$ or $R_{2^N}$;
ii) it can only cross the discontinuity manifolds at points on ${\BSigma}$, and 
iii) if it slides, it can only slide on  the intersection of all $N$ discontinuity manifolds, i.e., on ${\BSigma}$.
{\bf However}, the solution of a problem relative to perturbed initial conditions in general will not
satisfy the restricted motion described by points i)-iii) above, and it  may
slide on, or cross, some of the $\Sigma_i$'s and not just ${\BSigma}$. 
This fact must be taken into account when deriving the expression of the fundamental matrix
solution of the linearized dynamics, in particular of the saltation matrices.
\end{rem}

Next, we study the case $N=2$ in detail. The generalization to the case $N>2$ is simple, and
appropriate modifications required to describe the case $N>2$ are given below. 
Following the tree diagram for $N=2$, we have the following four subregions of  phase space
and corresponding vector fields, for ${\x}=\bmat{x_1\\x_2}\in \R^{2n}$:
\[
\begin{matrix} 
R_1=\{{\x}\,\ | \,\  
h_1({\x})<0, \,\ h_2({\x}) <0 \}, & F_1({\x})=\bmat{f^-(x_1) \\ f^-(x_2)}, \\
R_2=\{{\x} \,\ | \,\ h_1({\x})<0, \,\ h_2({\x}) >0 \}, & F_2({\x})=\bmat{f^-(x_1)\\ f^+(x_2)}, \\
R_3=\{{\x} \,\ | \,\ h_1({\x})>0, \,\ h_2({\x})<0 \}, & F_3({\x})=\bmat{f^+(x_1) \\ f^-(x_2)}, \\
R_4=\{{\x} \,\ | \,\ h_1({\x})>0, \,\ h_2({\x})>0 \}, & F_4({\x})=\bmat{f^+(x_1) \\ f^+(x_2)},
\end{matrix}
\]
and  $F({\x})$ in \eqref{PWSnetwork} is equal to $F_i({\x})$ for ${\x} \in R_i$.
Moreover we have $\Sigma_{1,2}=\{{\x} \in \R^{2n} \,\ | \,\ h_{1,2}({\x})=0\}$ and 
we consider also the sets $\Sigma_{1,2}^\pm$, defined as follows:
$\Sigma_1^\pm=\{{\x} \in \R^{2n} \,\ | \,\ h_{1}({\x})=0$, and $h_2(\x) \gtrless 0 \}$,
and similarly for  $\Sigma_2^\pm$.
The synchronous solution is ${\x_S}(t)=\bmat{x_S(t)\\ x_S(t)}$, and (under
Assumption \ref{xs_ass}) it evolves schematically as on the right of
Figure \ref{po_types}.

As already pointed out in Remark \ref{synch_rem}, while $\x_S(t)$ can only evolve in $R_1$ and $R_4$,
and can only cross/slide-on ${\BSigma}=\Sigma_1 \cap \Sigma_2$,  a perturbed solution might
instead cross just $\Sigma_1$ or $\Sigma_2$, and evolve in $R_2$ or $R_3$ 
and it might slide along $\Sigma_1$ and/or 
$\Sigma_2$. Hence, in order to compute the fundamental matrix solution, we will need the 
expressions of the sliding vector fields on $\Sigma_1$, $\Sigma_2$ and ${\BSigma}$.  
Let $f_{\Sigma}: \R^n \to \R^n$ denote the sliding vector field of \eqref{nonsmooth_eq} on 
$\Sigma$ (see \eqref{SlidingDyn}) and let $F_{{\mathcal S}}: \R^{2n} \to \R^{2n}$ denote 
the sliding vector field on  
${\mathcal S}$, where ${\mathcal S}$ is any of the following:
${\mathcal S}=\Sigma_1^\pm, \Sigma_2^\pm,{\BSigma}$. 

The next result provides  --at a point $\x$ on $\bf \Sigma$-- 
the four vector fields on $\Sigma_{1,2}^\pm$: these are the vector fields {\emph {that are felt}} at $\x$.

\begin{lemma}\label{OnSigma_prop}
Let $\x\in \bf \Sigma$, $\x=\bmat{x \\ x}$, and suppose that $\x$ is an attractive
point with respect to both $\Sigma_1$ and $\Sigma_2$: cfr \eqref{SlidePt}.  That is,
we have
\begin{equation}\label{PartAtt}
 0<\frac{\nabla h(x)^Tf^-(x)}{\nabla h(x)^T(f^-(x)-f^+(x))}<1 \ .
\end{equation}
Then, we have
\begin{equation*}\label{sliding_vf_eq}
F_{\Sigma_1^\pm}({\x})=\begin{pmatrix} f_{\Sigma}(x) \\ f^{\pm}(x) \end{pmatrix}, \,\
F_{\Sigma_2^\pm}({\x})=\begin{pmatrix} f^{\pm}(x) \\ f_{\Sigma}(x) \end{pmatrix},\,\
F_{{\BSigma}}({\x})=\begin{pmatrix} f_{\Sigma}(x) \\ f_{\Sigma}(x) \\ \end{pmatrix}=
e \otimes f_{\Sigma}(x), 
\end{equation*} 
where with $f_{\Sigma}(x)$ we denote the sliding vector field of \eqref{nonsmooth_eq} on
$\Sigma$, as defined in \eqref{SlidingDyn}.
\end{lemma}
\begin{proof}
We prove the statement for $F_{\Sigma_2^-}$. The proofs for the other sliding vector fields are analogous. 

Let $\x^k$ be a sequence of points in $\Sigma_2^-$, converging to $\x$: 
$\displaystyle{\lim_{k\to \infty}}\x^k =\x\in {\BSigma}$.  Since $\x^k\in \Sigma_2^-$,
we have $\x^k=\bmat{x_1^k \\ x_2^k}$ and $h(x_2^k)=0$ for all $k$.  

Now, letting 
$\hat F_1(\x^k)=F_1(\x^k)+\sigma M\x^k=
\bmat{f^-(x_1^k) \\ f^-(x_2^k)}+\sigma M\x^k$ and 
$\hat F_1(\x^k)=F_2(\x^k)+\sigma M\x^k=
\bmat{f^-(x_1^k) \\ f^+(x_2^k)}+\sigma M\x^k$, a Filippov
sliding vector field on $\Sigma_2^-$ is given by the convex combination
$F_{\Sigma_2^-}=(1-\alpha^k)\hat F_1(\x^k)+\alpha^k \hat F_2(\x^k)$, where 
$\alpha^k$ must be found from the requirement that $F_{\Sigma_2^-}$ is on
tangent plane, that is
$$(0,\nabla h^T)\left[ (1-\alpha^k)\bmat{f^-(x_1^k) \\ f^-(x_2^k)}+\alpha^k
\bmat{f^-(x_1^k) \\ f^+(x_2^k)}+\sigma M\x^k\right] \ = \ 0\ .$$
Therefore, 
$$\alpha^k=\frac{\nabla h^Tf^-(x_2^k)+\sigma(0,\nabla h^T)M\x^k}{\nabla h^T(f^-(x_2^k)-f^+(x_2^k))}\ .$$
Now, for any $\epsilon>0$,
there exists $K_\epsilon$ such that if $k>K_\epsilon$, then $\|\x^k-\x\|<\epsilon$.  This means that
(for $k$ sufficiently large) all points $\x^k$ are attractive sliding points relative to $\Sigma_2^-$, since
--because of \eqref{PartAtt}-- $0<\alpha^k<1$ for $k$ sufficiently large.
Thus, the sequence of sliding vector fields on $\Sigma_2^-$ is well defined. In the limit
as $k\to \infty$, $\alpha^k \to \frac{\nabla h^Tf^-(x)}{\nabla h^T(f^--f^+)(x)}$ 
(since $M\x=0$ for $\x\in \bf \Sigma$) and hence (see  \eqref{SlidingDyn}) we obtain 
$$F_{\Sigma_2^-}(\bf x) \ =  \ \bmat{f^-(x) \\ f_{\Sigma}(x)}.$$
\end{proof}

\begin{rem}
The result {\bf extends} to the case $N>2$ as follows.
The sliding vector field $F_{\Sigma_i}({\x})$ along a single 
$\Sigma_i$, has $f_{\Sigma}(x)$ in the $i$-th block while the components of the vector field 
in the $j$-th block, $j \neq i$, 
are equal to $f^+(x)$ or $f^-(x)$ in agreement with the sign of $h_j(x)$.
Moreover, $F_{\BSigma}({\x})= 
\begin{pmatrix}f_{\Sigma}(x) \\ \vdots \\ f_{\Sigma}(x) \end{pmatrix}$.
\end{rem}

Next, in the case of $N=2$, let $X(T,0)$ denote the monodromy matrix along $\x_S(t)$.
Then, because of Assumption \ref{xs_ass}, $X(T,0)$ can be written as: 
\begin{equation}\label{monodromy_eq}
X(T,0)=  X_1(T,t_3)X_{\BSigma}(t_3,t_2)S_{4,{\BSigma}}X_4(t_2,t_1)S_{14}X_1(t_1,0),
\end{equation}  
where
\[
\begin{matrix}
\frac{d X_i(t,\tau)}{dt}= & (DF_i({\x_S}(t))+\sigma M)X_i(t,\tau), & X_i(\tau,\tau)=I_{2n}, & i=1,2, \\
& & & \\
\frac{d X_{\BSigma}(t,\tau)}{dt}= & (DF_{\BSigma}({\x_S}(t))+\sigma M)X_{\BSigma}(t,\tau), 
& X_{\BSigma}(\tau,\tau)=I_{2n}. & \\
\end{matrix}
\]
In \eqref{monodromy_eq}, the matrix $S_{14}$ is the jump matrix from $R_1$ into $R_4$,
and the matrix $S_{4,{\BSigma}}$ 
is the jump matrix from $R_4$ into ${\BSigma}$. In general, see \cite{Ivanov} and \cite{Dieci.Lopez_2011},
the jump matrix on the  intersection of two discontinuity manifolds is ambiguous.  However,
in Theorems \ref{S14_thm}  and \ref{S4Sigma_thm} we show that for our problem
we can use a unique expression for the jump matrix.

\begin{rem}
The {\bf extension} of \eqref{monodromy_eq} {\bf for $N>2$}, is easily obtained
replacing $S_{14}$ and $S_{4,{\BSigma}}$ with $S_{1,2^N}$ and 
$S_{2^N,{\BSigma}}$.
\end{rem}

Now, let $\y(t)=\bmat{y_1(t) \\ y_2(t)}$ be the solution of \eqref{PWSnetwork} with 
perturbed initial
conditions $\y(0)={\s}_0+\Delta {\s}$, $\|\Delta \s\| \ll 1$.
Recall that the monodromy matrix expresses, at
first order in $\Delta {\s}$, the evolution of these perturbed initial conditions after one period.
Below we give details on how to compute $S_{14}$.   To go from $R_1$ to $R_4$, 
$\y(t)$ might either cross directly ${\BSigma}$ at time 
$t_1+\Delta t$ or it might instead cross $\Sigma_1$ and $\Sigma_2$ at two different times 
$\Delta t_1$ and $\Delta t_2$ before entering $R_4$ (see Remark \ref{synch_rem}). In \cite{Ivanov},
the form of $S_{14}$ is given and the author points out that the jump matrix is ambiguous.
In Lemma \ref{S14_prop} and 
\ref{S14_2_prop} we give the expression of $S_{14}$ for the two different possibilities when 
$\y(t)$ crosses $\Sigma_1$ and $\Sigma_2$ at two different
times, or it crosses directly ${\BSigma}$; but, then, Lemma \ref{S14_2_prop} states that only one expression of $S_{14}$ is
needed to assess perturbations of synchronous solutions and in Theorem \ref{S14_thm} we give
this unique expression of $S_{14}$ as a Kronecher product of the identity matrix $I_2$ with the jump 
matrix of the single agent.

\begin{lemma}\label{S14_prop} 
Let $\x_S=\bmat{x_S\\x_S}$ where $x_S$ satisfies Assumption \ref{xs_ass}.
Let	$\y(t)=\bmat{y_1(t)\\ y_2(t)}$ 
be the solution of \eqref{PWSnetwork} with initial
conditions $\y(0)={\s}_0+\Delta {\s}$, $\s_0=\bmat{s_0\\ s_0}$ and
$\|\Delta {\s}\|\ll 1$.
Let $(t_1+\Delta t_1)$ and $(t_1+\Delta t_2)$ be the times at which $\y(t)$ crosses respectively 
$\Sigma_1$ and $\Sigma_2$ before entering $R_4$.  Assume $\Delta t_1 \neq \Delta t_2$. 
Then, the jump matrix $S_{14}$ is given by 
\[S_{14}=I_2 \otimes I_n +I_2 \otimes \frac{(f^+(s_1)-f^-(s_1) )\nabla h(s_1)^T}{\nabla h(s_1)^Tf^-(s_1)}=I_2 \otimes S_{-,+},\]
with $S_{-,+}=I_n+\frac{(f^+(s_1)-f^-(s_1) )\nabla h(s_1)^T}{\nabla h(s_1)^Tf^-(s_1)}$ being the jump matrix for the single agent in \eqref{nonsmooth_eq} from $R^-$ to $R^+$.
\end{lemma}
\begin{proof}
We examine below the case $\Delta t_1<\Delta t_2$. The proof for the case $\Delta t_2>\Delta t_1$ is similar. 
Following \cite{Ivanov} and \cite{Dieci.Lopez_2011}, we can rewrite $S_{14}$ as the composition of two 
jump matrices: $S_{14}=S_{34}S_{13}$,  with $S_{ij}$ being the jump 
matrix from $R_i$ to $R_j$.  Using standard 
results on jump matrices for one discontinuity manifold, we have that 
(note that $\sigma M\x=0$ along a synchronous solution)
\[S_{13}= \begin{pmatrix} S_{-,+} & 0  \\ 0 & I_n \end{pmatrix} ,\,\
S_{34}= \begin{pmatrix} I_n & 0 \\ 0 & S_{-,+} \end{pmatrix},\]
with $S_{-,+}$ defined in the statement. Then $S_{34}S_{13}=I_2 \otimes S_{-,+}.$ 
The case $\Delta t_1>\Delta t_2$, in which the perturbed solution first crosses $\Sigma_2$ to enter into 
$R_2$ and then crosses $\Sigma_1$ to enter into $R_4$, gives $S_{14}$ as the product $S_{24}S_{12}$. 
It is easy to verify that 
$S_{24}S_{12}=S_{34}S_{13}$. 
\end{proof}

In the Lemma below we consider the case when the perturbed solution $\y(t)$ crosses $\Sigma_1$ and $\Sigma_2$ at the same time, and 
give a unique expression for the jump matrix $S_{14}$ in this case.

\begin{lemma}\label{S14_2_prop}	
As in Lemma \ref{S14_prop}, let $\x_S=\bmat{x_S\\x_S}$ where $x_S$ satisfies Assumption \ref{xs_ass}.
Let	$\y(t)=\bmat{y_1(t)\\ y_2(t)}$ 
be the solution of \eqref{PWSnetwork} with initial
conditions $\y(0)={\s}_0+\Delta {\s}$, $\s_0=\bmat{s_0\\ s_0}$ and
$\|\Delta {\s}\|\ll 1$.
Let $(t_1+\Delta t_1)$ and $(t_1+\Delta t_2)$ be the times at which $\y(t)$ crosses respectively 
$\Sigma_1$ and $\Sigma_2$ before entering $R_4$.  
Assume $\Delta t_1=\Delta t_2=\Delta t$, so that
necessarily ${\y}(t)$ from $R_1$ crosses $\BSigma$ to enter into $R_4$.
Then, without loss of generality,   
\[S_{14}=\begin{pmatrix} I_n+\frac{(f^+(s_1)-f^-(s_1)) \nabla h^T(s_1)}{\nabla h^T(s_1)^Tf^-(s_1)} & {\bf 0} \\		
	\frac{(f^+(s_1)-f^-(s_1)) \nabla h^T(s_1)}{\nabla h^T(s_1)f^-(s_1)} & I_n		
\end{pmatrix}=I_2 \otimes I_n+\begin{pmatrix} 1 &0 \\ 1 & 0 \end{pmatrix} \otimes 
\frac{(f^+(s_1)-f^-(s_1))\nabla h^T(s_1)}{\nabla h^T(s_1)f^-(s_1)}.\]	
\end{lemma}
\begin{proof}
Following  \cite{Ivanov} and \cite{Dieci.Lopez_2011}, 
there are two possible expressions for $S_{14}$ when $\Delta t_1=\Delta t_2$: 
\begin{eqnarray*}
	S_{14}^{(1)}=I_{2n}+\left (F_4({\bf s}_1)-F_1({\bf s}_1)\right ) 
	\frac{\nabla h_1^T(s_1)}{\nabla h_1(s_1)^TF_1({\bf s}_1)}=
	\begin{pmatrix} I_n+\frac{(f^+-f^-) \nabla h^T(s1)}{\nabla h^Tf^-(s_1)} & {\bf 0} \\
		\frac{(f^+-f^-)  \nabla h^T(s1)}{\nabla h^Tf^-(s_1)} & I_n
	\end{pmatrix} ,\\
	S_{14}^{(2)}=I_{2n}+\left (F_4({\s}_1)-F_1({\s}_1) \right )
	\frac{\nabla h_2^T(s_1)}{\nabla h_2(s_1)^TF_1({\bf s}_1)}=
	\begin{pmatrix} I_n & \frac{(f^+-f^-) \nabla h(s_1)2^T(s1)}{\nabla h(s_1)^Tf^-(s_1)} \\
		{\bf 0 } & I_n+\frac{(f^+-f^-) \nabla h^T(s1)}{\nabla h^T(s_1)f^-(s_1)}  
	\end{pmatrix}.
\end{eqnarray*} 
Below we show that the two expressions are equivalent, in the sense that their
action on the input vector $\x_S(t_1)-\y(t_1)$ is identical.

Explicit computation of $S^{(1)}_{14} ({\x_S}(t_1)-{\y}(t_1))$ and 
$S^{(2)}_{14}({\x_S}(t_1)-{\y}(t_1))$ gives  
\begin{align*}
S_{14}^{(1)} ({\x_S}(t_1)-{\y}(t_1))= &
\begin{pmatrix} 
	(s_1-y_1(t_1))+(f^+(s_1)-f^-(s_1))
	\frac{ \nabla h(s_1)^T(s_1-y_1(t_1))}{\nabla h(s_1)^Tf^-(s_1)} \\
	(f^+(s_1)-f^-(s_1))\frac{ \nabla h(s_1)^T(s_1-y_1(t_1))}{\nabla h(s_1)^Tf^-(s_1)}+(s_1-y_2(t_1))
\end{pmatrix},\\
S_{14}^{(2)} ({\x_S}(t_1)-{\y}(t_1))= &
\begin{pmatrix} 
	(s_1-y_1(t_1))+(f^+(s_1)-f^-(s_1))\frac{\nabla h(s_1)^T(s_1-y_2(t_1))}{\nabla h(s_1)^Tf^-(s_1)} \\
	(f^+(s_1)-f^-(s_1))\frac{ \nabla h(s_1)^T(s_1-y_2(t_1))}{\nabla h(s_1)^Tf^-(s_1)}+(s_1-y_2(t_1))
\end{pmatrix}.
\end{align*}
The statement of the theorem then follows if we show 
\begin{equation}\label{S14(1)(2)_eq}
\nabla h(s_1)^T(s_1-y_1(t_1))=\nabla h(s_1)^T(s_1-y_2(t_1))+\text{h.o.t.}.
\end{equation}
We use the following Taylor expansions 
\begin{align*}
	y_{1,2}(t_1+\Delta t)& =y_{1,2}(t_1)+f^-(y_{1,2}(t_1))\Delta t_1+O(\Delta t_1^2), \\
	h(y_{1,2}(t_1))& =h(s_1)+\nabla h^T(s_1)(y_{1,2}(t_1)-s_1)+O(\|y_{1,2}(t_1)-s_1\|^2), \\
	0=h(y_{1,2}(t_1+\Delta t))&=h(y_{1,2}(t_1))+
	\nabla h(y_{1,2}(t_1))^Tf^-(y_{1,2}(t_1))\Delta t+\text{h.o.t.}=\\
	& =h(s_1)+\nabla h(s_1)^T(y_{1,2}(t_1)-s_1)+ \nabla h(s_1)^Tf^-(s1)\Delta t+\text{h.o.t.}.
\end{align*}
From the last equality \eqref{S14(1)(2)_eq} follows.
\end{proof}

We now have the following theorem, that gives a unique expression for $S_{14}$.
\begin{theorem}\label{S14_thm}
With same notation as in Lemma \ref{S14_2_prop}, 
without loss of generality we can use the following expression for the jump matrix $S_{14}$:  
\[S_{14}=I_2 \otimes S_{-,+},
\] 
where $S_{-,+}$ is the jump matrix of the single agent  \eqref{nonsmooth_eq} from $R^-$ into $R^+$. 
\end{theorem}
\begin{proof}
As a consequence of Lemmas \ref{S14_prop} and \ref{S14_2_prop}, the jump matrix 
$S_{14}$ can be taken as follows
\[S_{14}= \left \{ \begin{matrix} I_2 \otimes I_n +I_2 \otimes \frac{(f^+(s_1)-f^-(s_1) )\nabla h(s_1)^T}{\nabla h(s_1)^Tf^-(s_1)}, & \Delta t_1 \neq \Delta t_2, \\ 
I_2 \otimes I_n+\begin{pmatrix} 1 & 0 \\ 1 & 0\end{pmatrix} \otimes 
\frac{(f^+(s_1)-f^-(s_1)) \nabla h(s_1)^T}{\nabla h(s_1)^Tf^-(s_1)} ,  
 & \Delta t_1 = \Delta t_2. \end{matrix} \right .\] 
Moreover, if $\Delta t_1=\Delta t_2$, then \eqref{S14(1)(2)_eq} implies that 
at first order $\nabla h(s_1)^Ty_1(t_1)=\nabla h(s_1)^Ty_2(t_1)$ and similarly to the proof of Lemma \ref{S14_2_prop}, easy computations imply that 
 \[\left ( I_2 \otimes \frac{(f^+(s_1)-f^-(s_1)) \nabla h(s_1)^T}{\nabla h(s_1)^Tf^-(s_1)} \right ) \y(t_1)=
 \left ( \begin{pmatrix} 1 & 0 \\ 1 & 0 \end{pmatrix} \otimes 
 \frac{(f^+(s_1)-f^-(s_1)) \nabla h(s_1)^T}{\nabla h(s_1)^Tf^-(s_1)} \right ) \y(t_1).\] Then, at first order,
 \[\y(t_1^+)-{\x}(t_1^+) =(I_2 \otimes S_{-,+})(\y(t_1^-)-{\x}(t_1^-)), \]
 for any $\y(t)$ perturbed solution of ${\x_S}(t)$.
\end{proof}

\begin{rem}
Lemma \ref{S14_prop} and \ref{S14_2_prop} and Theorem \ref{S14_thm} {\bf extend}
in a straightforward 
way to the case $N>2$, via the replacement of $S_{14}$ with 
$S_{1,2^N}=I_N \otimes S_{-,+}$,  
giving $\y(t_1^+)-{\x}(t_1^+)=(I_N \otimes S_{-,+})(\y(t_1^-)-{\x}(t_1^-))$.
\end{rem}

Going back to \eqref{monodromy_eq}, we next need to analyze
$S_{4,{\BSigma}}$, that is the jump matrix from $R_4$ to $\bf \Sigma$.   In general, the 
jump matrix from the region $R_4$ to $\bf \Sigma$ is not uniquely defined, which makes it
not possible to give a unique expression for the mondromy matrix.   However, 
the jump matrix for ${\x_S}(t)$ is an exception, as stated in Theorem \ref{S4Sigma_thm},
in the sense that (at first order) we can give a unique expression for the action of
$S_{4,{\BSigma}}$ on $\y(t_2^-)-\x_S(t_2^-)$. The proof in Theorem \ref{S4Sigma_thm} sums up 
results analogous to the ones given in Lemmas \ref{S14_prop} and \ref{S14_2_prop} and in Theorem \ref{S14_thm} for $S_{14}$.

\begin{theorem}\label{S4Sigma_thm}
The jump matrix $S_{4,{\BSigma}}$ can be taken to be 
\[S_{4,{\BSigma}}=I_2 \otimes (I_n+\frac{(f_{\Sigma}(s_2)-f^+(s_2)) \nabla h(s_2)^T}{\nabla h(s_2)^Tf^+(s_2)})=I_2 \otimes S_{+,\Sigma},\] 
where $S_{+,\Sigma}$ is the jump matrix of the single agent \eqref{nonsmooth_eq} from $R^+$ into 
$\Sigma$, and it is explicitly given by 
$S_{+,\Sigma}=\displaystyle{I+\frac{(f_\Sigma(s_2)-f^+(s_2))\nabla h(s_2)^T}{\nabla h(s_2)^Tf^+(s_2)}}$, with 
$f_{\Sigma}$ given in \eqref{SlidingDyn}.
\end{theorem}  
\begin{proof}
Following \cite{Dieci.Lopez_2011}, a perturbed solution of ${\x_S}(t)$ 
might either reach directly ${\bf \Sigma}$ in a neighborhood of ${\s}_2$  or it might first slide along 
$\Sigma_1$ or $\Sigma_2$ before reaching ${\bf \Sigma}$.  Let 
$\y(t)=\bmat{y_1(t) \\ y_2(t)}$ be the perturbed solution, with $\y(0)={\s}_0+ \Delta {\s}$ and 
$\|\Delta {\s}\| \|\ll1$.
\begin{itemize}
\item[i)] We first consider the case in which the perturbed solution reaches $\Sigma_1^+$, slides along 
it and then it reaches ${\bf \Sigma}$. Then $S_{4, {\bf \Sigma}}=S_{{\bf \Sigma},\Sigma_1^+}S_{\Sigma_1^+,4}$ with 
\begin{align*}
S_{4,\Sigma_1^+} & =I_{2n}+\frac{(F_{\Sigma_1^+}({\s}_2)-F_4({\s}_2))\nabla h_1({\s}_2)^T}{\nabla h_1({\s}_2)^TF_1({\s}_2)}=
\begin{pmatrix} I_n+\frac{(f_{\Sigma}({s}_2) -f^+({s}_2))\nabla h(s_2)^T}{\nabla h(s_2)^Tf^+({s}_2)} & 0 \\
	0 & I_n  \end{pmatrix},  \\
S_{\Sigma_1^+,{\bf \Sigma}} & =I_{2n}+\frac{(F_{\Sigma}({\s}_2)-F_{\Sigma_1^+}({\s}_2))\nabla h_2({\s}_2)^T}{\nabla h_2({\s}_2)^TF_1({\s}_2)}=
\begin{pmatrix} I_n & 0 \\
	0 & I_n+\frac{(f_{\Sigma}({s}_2) -f^+(x))\nabla h(s_2)^T}{ \nabla h(s_2)^Tf^+({s}_2)} \end{pmatrix}, 
\end{align*}
where the vector fields $F_{{\bf \Sigma}}$ and $F_{\Sigma_1^+}$ are as in Lemma 
\ref{OnSigma_prop}. Then the statement follows in this case. The case in which $\y(t)$ 
reaches $\Sigma_2^+$ before reaching $\bf \Sigma$ is analogous and gives the same expression for the jump matrix.   

\item[ii)] We next consider the case in which $\y(t)$ reaches $\bf \Sigma$ directly at time 
$t_2+\Delta t$. As in \cite{Dieci.Lopez_2011}, in this case 
$S_{4,{\BSigma}}$ has two expressions, namely 
\begin{align*}
S_{4,{\BSigma}}^{(1)} & =I_{2n}+
\frac{(F_{{\BSigma}}-F_4)({\s}_2)\nabla h_1({\s}_2)^T}{\nabla h_1({\s}_2)^TF_4({\s})}, \\
S_{4,{\BSigma}}^{(2)} & =I_{2n}+
\frac{(F_{{\BSigma}}-F_4)({\s}_2)\nabla h_2({\s}_2)^T}{\nabla h_2({\s}_2)^TF_4({\s})}. 
\end{align*} 
Using the same reasonings as in Lemma \ref{S14_2_prop}, at first order
the perturbed solution must 
satisfy $\nabla h_1({\s}_2)^T\y(t_2)=\nabla h_2({\s}_2)^T\y(t_2)$, i.e., 
$\nabla h(s_2)^Ty_1(t_2)=\nabla h(s_2)^Ty_2(t_2)+\text{h.o.t.}$.
In this case we can use a unique expression for $S_{4,\BSigma}$:
\[S_{4,{\BSigma}}=I_2 \otimes I_n+\begin{pmatrix} 1 & 0 \\ 1 & 0 \end{pmatrix}
\frac{(f_{\Sigma}-f^+)(s_2)\nabla h^T}{\nabla h^Tf^+(s_2)}.\]
\item[iii)] The final argument is analogous to the one in the proof of Theorem \ref{S14_thm}. 
We have two possible expressions for the jump matrix:  in  i), when $\y$ does not reach $\bf \Sigma$ 
directly and in ii) when $\y$ reaches $\bf \Sigma$ directly from $R_4$. In this last case 
however, $\nabla h^T(s_2)y_1(t_2)=\nabla h^T(s_2)y_2(t_2)+\text{h.o.t.}$ and this implies that at first order 
\[\y(t_2^+)-\x(t_2^+)=(I_2 \otimes S_{+,\Sigma})(\y(t_2)-\x(t_2) ). \] 
\end{itemize}
The above points i), ii), iii) imply the statement of the theorem.
\end{proof}

\begin{rem}
The {\bf extension} of Theorem \ref{S4Sigma_thm} 
to the case of $N>2$ is immediately achieved by taking
$S_{2^N,{\BSigma}}=I_N \otimes S_{+,\Sigma}$.
\end{rem}  

Finally, the following theorem gives the complete
expression for the monodromy matrix of \eqref{PWSnetwork} linearized about the synchronous
solution $\x_S(t)$.  The proof puts together all the results previously derived in this section
and is therefore omitted.
\begin{theorem}\label{X_thm}
Let $x_s(t)$ be limit cycle of \eqref{nonsmooth_eq}, and let Assumption \ref{xs_ass} be satisfied.
Let ${\x_S}(t)=e \otimes x_S(t)$ be the corresponding synchronous solution of \eqref{PWSnetwork}.
Then, the monodromy matrix of \eqref{PWSnetwork} along $\x_S$ can be taken as
\[X(T,0)=X_1(T,t_3)X_\Sigma(t_3,t_2)(I_N \otimes S_{+,\Sigma})X_{2^N}(t_2,t_1)(I_N \otimes 
S_{-,+})X_1(t_1,0),
\] 
with 
\[
\begin{matrix}
\frac{d X_i(t,\tau)}{dt}= & (DF_i({\x_S}(t))+\sigma M)X_i(t,\tau), & X_i(\tau,\tau)=I_{Nn}, & i=1,2^N \\
& & & \\
\frac{d X_{{\BSigma}}(t,\tau)}{dt}= & (DF_{{\BSigma}}({\x_S}(t))+\sigma M)X_{{\BSigma}}(t,\tau), & 
X_{{\BSigma}}(\tau,\tau)=I_{Nn}, & \\
\end{matrix}
\]
where $DF_i$  and $DF_{{\BSigma}}$ are the Jacobian matrices of the vector fields 
$F_i$, $i=1,2^N$ and $F_{{\BSigma}}$ respectively. 
\hfill\qed
\end{theorem}

The key implication of Theorem \ref{X_thm} is that 
the saltation matrices appearing in the expression of $X(T,0)$ can be obtained from the saltation 
matrices of a single agent in the network, which is a great simplification.  However,  computation of the 
matrices $X_i$ and $X_{{\BSigma}}$ involves all $N$ agents. 
In case of large networks, these computations are too expensive,
and an extension of the MSF theory to \eqref{PWSnetwork} is needed.    This is the
purpose of the next section.

\section{Master Stability function}\label{MSF}
The MSF is a very nice technique which allows to study linearized stability of the synchronous solution
of a smooth network of $N$ agents of size $n$ each, by working with $N$ linearized systems
of size $n$, rather than one linearized system of size $nN$, a substantial saving!
The key idea in the MSF technique is to consider the variational equation along the synchronous
solution and to perform a change of coordinates induced by the matrix of eigenvectors of the
Laplacian $L$.   For smooth systems,
this change of coordinates brings the whole network into a block diagonal structure with sub-blocks 
of size $n \times n$  and this in turn allows one to study the stability of $N$ systems of dimension
$n$ instead of the stability of one system of dimension $nN$ (see \cite{Pecora.Carroll}). 

In the case of piecewise smooth vector fields for the agents, there are
at least two new concerns.  First, the same change of coordinates, while still
bringing the  variational equation along the synchronous solution into block diagonal form, will also change
the equations of the discontinuity manifolds that will now in general involve more than one agent,
and possibly all of them.  Then, we should not expect the saltation matrices to be block 
diagonal. However,  we will see
that the transformation preserves the block structure of the 
saltation matrices of Theorem \ref{X_thm}.  In other words, the Kronecker products involved in the expression of the
saltation matrices along a synchronous solution 
are left unchanged by the coordinate change induced by the eigenvectors
of $L$ (see equation \eqref{jump_unchanged_eq}).  
The second concern is related to the portion of the fundamental matrix
solution on $\bf \Sigma$, in particular to the Jacobian of $F_{\BSigma}$.  
The general lack of uniqueness in expressing the sliding vector
field on the intersection of two or more discontinuity manifolds is not a concern in this setting, 
since, by Lemma \ref{OnSigma_prop},  on
$\bf \Sigma$ we have a unique Flilippov sliding vector field.  However a difficulty is
related to expressing the Jacobian itself, since again the coordinate change 
seemingly will destroy the sought block structure.  We deal with this difficulty in 
Lemma \ref{Kron_slide_prop} below. 
As a final result, in Theorem  \ref{main_thm} we will see that we can use the MSF technique also for PWS networks, and in particular to study stability of the synchronous solution $\x_S$ by linearized analysis on $N$ systems of size $n$.

In what follows, let $L$ be the graph Laplacian matrix ($L=L^T$),
and $W$ be the matrix of the orthonormal eigenvectors of $L$:  
$W^TLW=\Lambda$, with $\Lambda$ diagonal. 
Then $M=L \otimes E$ has eigenvalues $ \lambda_i\mu_j$, $i=1, \ldots N$, $j=1, \ldots, n$,  
with $\lambda_{i}$'s the eigenvalues of $L$, and $\mu_j$'s the eigenvalues of $E$. 
Moreover, with $V=W \otimes I$,  
we get $V^{-1}MV=(W^T \otimes I_n)(L \otimes E)(W \otimes I_n)= \Lambda \otimes E$. 
In \cite{Pecora.Carroll}, the change of variables $\y=V^{-1}{\x}$ 
reduces the variational equation along a synchronous solution
$\x_S(t)$ into the following block diagonal form
\[\dot z_i=(Df(x_S)+\sigma \lambda_{i}E) z_i, \qquad i=1, \ldots, N,\]
that is we have $N$ systems of size $n$.
The issue with nonsmooth agents is that the change of variables above in general changes the equations of 
the discontinuity manifolds as well, and this makes it impossible to study $N$ systems independently.
In particular, the new equations of the discontinuity manifold(s) might involve all the agents.  
\begin{exam}
To illustrate the last statement, take $L=\begin{pmatrix} -1 & 1 \\ 1 & -1 \end{pmatrix}$ and 
$M=L \otimes I_2$. Then, using the same notations as before, 
$W=\frac{1}{\sqrt{2}}\begin{pmatrix} 1 & 1 \\ 1 & -1 \end{pmatrix}$ and 
$V=W \otimes I_2$.  If the discontinuity surface for the single agent is the plane 
$\{ x \in \R^n | h(x)=e_1^Tx=0 \}$ , then $\Sigma_1$ and $\Sigma_2$ in the 
$\y$ coordinates become respectively 
$e_1^T\y=-e_3^T\y$ and $e_1^T\y=e_3^T\y$ .
It follows that the agents cannot be studied independently even though the variational 
equations for the fundamental matrix solution are in block diagonal form. 
\end{exam}

Nevertheless, we will now see that the special structure of the fundamental matrix solution for the case of 
a synchronous periodic solution of \eqref{PWSnetwork}, allows to study the stability
of $\x_S$ via $N$ systems of dimension $n$.  
   
With $V=W \otimes I$, consider the monodromy matrix for the linearization along the synchronous 
periodic solution.  Let $\y_S(t)=V^{-1}\x_S(t)$ and
$Y(T,0)=V^{-1}X(T,0)V$. Then 
\begin{align*}
Y(T,0) =Y_1(T,t_3)Y_{\BSigma}(t_3,t_2)(V^{-1}S_{2^N,{\BSigma}}V) 
Y_{2^N}(t_2,t_1)(V^{-1}S_{1,2^N}V)Y_1(t_1,0), 
\end{align*}
with $Y_i(t,\tau)=V^{-1}X_i(t,\tau)V$, $i=1,2^N$, and 
$Y_{{\BSigma}}(t,\tau)=V^{-1}X_{{\BSigma}}(t,\tau)V$. 

First notice that the particular structure of the jump matrices is such that 
\begin{eqnarray}\label{jump_unchanged_eq}
(V^{-1}S_{1,2^N}V)= & (W^T \otimes I_n) (I_N \otimes  S_{-,+})(W \otimes I_n)= \\
&  (W^T I_N W) \otimes ( I_n S_{-,+} I_n) =I_N \otimes S_{-,+},
\end{eqnarray}
and similarly for $S_{2^N,{\bf  \Sigma}}$.

Secondly, we show that the $Y_i$'s can be obtained solving block diagonal 
systems of ODEs.  
Note that 
\begin{align*}
V^{-1} DF_i(V{\x_S}(t))V & =(W^{T} \otimes I_n)(I_N \otimes Df_*(x_S(t))(W \otimes I_n) \\
& =I_N \otimes Df_*(x_S(t)), \,\ 
\end{align*}
with $Df_*=Df^-$ for $i=1$, and $Df_*=Df^+$ for $i=2^N$, and hence
\begin{equation}\label{Yi_eq1}
\frac{d Y_i(t,\tau)}{dt}= \left[I_N \otimes Df_*(x_S(t))+\sigma \Lambda  \otimes E\right]Y_i(t,\tau),   \quad Y_i(\tau,\tau)=I_{Nn},    \quad  i=1,2^N.
\end{equation}
The following Lemma shows how to rewrite the sliding vector field 
$V^{-1}DF_{\BSigma}(\x_S)V$ as a Kronecker product as well. 
\begin{lemma}\label{Kron_slide_prop}
	Let $\x_S(t)$ be a synchronous periodic solution of \eqref{PWSnetwork} and let 
	$F_{\BSigma}(\x_S(t))$ be the sliding vector field defined as in Lemma 
	\ref{OnSigma_prop}. Then 
\begin{equation}\label{DerSlideFMS}\begin{split}
V^{-1} DF_{\BSigma}(\x_S(t))V & =I_N \otimes Df_{\Sigma}(x_S(t))+  \\
&	\{\frac\sigma{\nabla h^T(f^--f^+)} \Lambda \otimes [(f^+-f^-)\nabla h^TE]\}_{x_S(t)}
+\sigma M\ . 
\end{split}\end{equation}
\end{lemma}
\begin{proof}
For simplicity, we will show the statement for $N=2$, the generalization for $N>2$ is immediate.  
	The sliding vector field on $\Sigma$ for \eqref{PWSnetwork} is defined in
	Lemma \ref{OnSigma_prop}.	That is, $F_{\BSigma}$ is 
\begin{equation*}
		F_{\BSigma}(\x) = \left [ \begin{matrix} 
			(1-\alpha_1(\x)))f^-(x_1)+\alpha_1(\x)f^+(x_1) \\
			(1-\alpha_2(\x)))f^-(x_2)+\alpha_2(\x)f^+(x_2)
		\end{matrix}\right ] +\sigma M \x_S,
\end{equation*}
where we have kept the term $M\x_S$ even though it is $0$ (since $\x_S$
is synchronous), to clarify the computation of the Jacobian.

Now, $\alpha_1({\x_S})$ and $\alpha_2(\x_S)$ must be chosen so that: 
$\bmat{\nabla h \\ 0}^TF_{\BSigma}({\x_S})=\bmat{0 \\ \nabla h}^TF_{\BSigma}({\x_S})=0$. 
Let $\alpha(x_S)$ be such that $f_{\Sigma}(x_S)=(1-\alpha(x_S))f^-(x_S)+\alpha(x_S)f^+(x_S)$ as in equation \eqref{SlidingDyn}. 	
Then, for $i=1,2$, when we impose the tangency conditions, using $M=L \otimes E$, we get:
	\[\alpha_i({\x_S})=\alpha(x_S)+\frac\sigma{\nabla h(x_S)^T(f^--f^+)(x_S)}\nabla h(x_S)^T(l_{i1}Ex_1+l_{i2}Ex_2)_{x_{1,2}=x_S}=\alpha(x_S),\]
where the last equality follows from the definition of the Laplacian matrix ($l_{i1}=-l_{i2}$ ) and the
fact that the solution is synchronous, i.e. $x_1=x_2=x_S$. 
	The gradient of $\alpha_1(\x)$ is then obtained as follows 
\begin{equation*}\begin{split}
D_{x_1}\alpha_1(\x_S)&=\left[\nabla \alpha(x_1)\right]_{x_1=x_S}	
+\frac\sigma{\nabla h^T(f^--f^+)(x_S)}l_{11}\nabla h(x_S)^TE+ \\
& \,\ \frac{\nabla h^T(l_{11}Ex_1+l_{12}Ex_2 )(\nabla h^T(Df^+-Df^-))}
{(\nabla h^T(f^--f^+))^2} \vert_{x_1=x_s}\\
		&= \left[\nabla \alpha(x_1)\right]_{x_1=x_S}
		+\frac\sigma{\nabla h^T(f^--f^+)(x_S)}l_{11}\nabla h(x_S)^TE ,
\end{split}\end{equation*}
and similarly
\begin{equation*}
D_{x_2}\alpha_1(\x_S)=\frac\sigma{\nabla h(x_S)^T(f^--f^+)(x_S)}l_{12}\nabla h(x_S)^TE .
\end{equation*}
Therefore,
$$
\nabla \alpha_1({\x_S})=\bmat{ \left[\nabla \alpha(x_1)\right]_{x_1=x_S}
+\frac\sigma{\nabla h^T(f^--f^+)}l_{11}\nabla h^TE \\
 \frac\sigma{\nabla h^T(f^--f^+)(x_S)}l_{12}\nabla h^TE}_{x_S}\ . 
$$
For $\alpha_2({\x_S})$ we obtain in a similar way
$$
\nabla \alpha_2({\x_S})=\bmat{\frac\sigma{\nabla h^T(f^--f^+)}l_{21}\nabla h^TE \\
\left[\nabla \alpha(x_2)\right]_{x_2=x_S}+\frac\sigma{\nabla h^T(f^--f^+)(x_S)}l_{22}\nabla h^TE}_{x_S}\ .  
$$
Then 
\begin{equation}\begin{split}\label{JacOnSigma}
		DF_{\Sigma}({\x_S}) & =\left [ \begin{matrix}
			Df_{\Sigma}(x_S)+\sigma l_{11}B &  \sigma l_{12}B \\
			& \\
			\sigma l_{21}B &  Df_{\Sigma}(x_S)+\sigma l_{22}B
		\end{matrix}
		\right ] +\sigma M \\
		& = I_N \otimes Df_{\Sigma}(x)+ \sigma L \otimes B+\sigma L\otimes E,
\end{split}\end{equation}
	where $B=\frac{(f^+-f^-)(x_S)\nabla h(x_S)^T}{\nabla h(x_S)^T(f^--f^+)(x_S)}E$. Then the statement follows at once.  
\end{proof}
\begin{rem}
In the case of $N=2$, the Laplacian matrix is trivial: $L=\bmat{-1 & 1 \\ 1 & -1}$.  But. aside
from this simplification, 
the proof of Lemma \ref{Kron_slide_prop} is identical for the case of $N>2$ except that we
have to account for more entries in the Laplacian matrix.   In particular, the
expression \eqref{JacOnSigma} remains valid.
\end{rem}

From Lemma  \ref{Kron_slide_prop} and \eqref{Yi_eq1},
with $Df_*=Df^-$ for $i=1$, and $Df_*=Df^+$ for $i=2^N$, we get
\begin{equation}\label{Yi_eq}\begin{split}
	 \frac{d Y_i(t,\tau)}{dt}& = \bigl(I_N \otimes Df_*(x_S(t))+\sigma \Lambda \otimes E\bigr)Y_i(t,\tau), \\ & Y_i(\tau,\tau)=I_{Nn},  \,\ i=1,2^N, \\
	\frac{d Y_{\Sigma}(t,\tau)}{dt}&= \bigl(I_N \otimes Df_{\Sigma}(x_S(t))+\sigma \Lambda \otimes E
	+\sigma \Lambda \otimes B	\bigr)Y_{\Sigma}(t,\tau),  \\ &  Y_{\Sigma}(\tau,\tau)=I_{Nn}. 
\end{split}\end{equation}
To sum up
\begin{equation}\label{Y_eq}
Y(T,0)=Y_1(T,t_3)Y_{\bf  \Sigma}(t_3,t_2)(I_N \otimes S_{+,{\bf  \Sigma}})
Y_{2^N}(t_2,t_1)(I_N \otimes S_{-,+})Y_1(t_1,0),
\end{equation}
with $Y_1$, $Y_{2^N}$ and $Y_{\BSigma}$ as in the block diagonal equations \eqref{Yi_eq}. 
In conclusion, we proved the following key theorem.
\begin{theorem}\label{main_thm}
The Floquet multipliers of \eqref{PWSnetwork} along the synchronous periodic orbit 
${\x_S}(t)$ are the $nN$ eigenvalues of the $N$ matrices $Z_i(T)$, 
$i=1, \ldots, N$, satisfying the following variational equations  
\begin{align}\label{Zi_eq}
\dot Z_i=  \left  \{ \begin{matrix}  
(Df^-(x_S(t))+\nu E) Z_i, & 0 \leq t<t_1, \\
(Df^+(x_S(t))+\nu E)Z_i,  & t_1<t<t_2, \\
(Df_{\Sigma}(x_S(t))+\nu  [E+B]) Z_i, &  t_2<t<t_3, \\
(Df^-(x_S(t))+\nu E) Z_i,&  t_3<t<T, 
\end{matrix} \right .
\end{align}
where $\nu=\sigma\lambda_i$, $i=1,\dots, N$, 
$B=\frac{(f^+-f^-)(x_S)\nabla h(x_S)^T}{\nabla h(x_S)^T(f^--f^+)(x_S)}E$, and 
subject to the initial conditions: $Z_i(0)=I_n$,  
$Z_i(t_1^+)=S_{-,+}Z_i(t_1^-)$,  $Z_i(t_2^+)=S_{+,{\BSigma}}Z_i(t_2^-)$.
\hfill\qed
\end{theorem}
As usual, we call Floquet exponents the logarithms of the multipliers.  For sure there is a $0$
exponent, since $1$ is a multiplier, because $\x_S$ is a periodic solution of \eqref{PWSnetwork}.
Now, since our network is connected, $L$ has only one 
eigenvalue equal to $0$, let it be $\lambda_1=0$, all other eigenvalues of $L$ being
negative.  With this observation we are ready for the following definition.

\begin{definition}\label{msf_def}
Let $\tau_{i,j}(\sigma)$ be the multipliers of \eqref{Zi_eq}, for $\nu=\sigma \lambda_i$, and
$i=2,\dots, N$, and let $l_{i,j}=\log|\tau_{i,j}|$, $j=1,\dots, n$.
The Master Stability Function (MSF)
for \eqref{PWSnetwork}, relative to the synchronous periodic solution $\x_S$, is the
largest value $l_{i,j}$, call it $\lambda$.  The synchronous manifold
(i.e., the synchronous solution) $\x_S$ is {\sl transversally stable} for those values
of $\sigma$, if any, for which $\lambda<0$.
\end{definition}

\begin{rem}
In the literature for smooth networks, the MSF is defined in terms of the Lyapunov
exponents of the linearized problem.  Of course, in the case of periodic orbits, these
are the Floquet exponents, and hence our definition is consistent with previous
usage of the MSF.
\end{rem}

Naturally, the value of the MSF $\lambda$ depends on the coupling strength $\sigma$,
as well as on $E$ and the Laplacian $L$.  However, for a given network topology (hence, for
given $L$ and $E$), the MSF depends
only on $\sigma$.   We must further appreciate that the network is synchronizable about 
$\x_S$, if all parameters values $\sigma\lambda_k$, $k=2,\dots, N$, give multipliers
less than $1$ in modulus.
Moreover, this is true regardless of whether or not $x_S$ is an asymptotically
stable periodic orbit of the system \eqref{nonsmooth_eq}.

\section{Periodic orbit of a piecewise smooth mechanical system. \\ Computation of 
the MSF}\label{GalvaSection}
Here we study a system of two identical piecewise smooth mechanical oscillators,
first studied in \cite{Galvanetto_95}.
When they are not coupled, the single agents have an asymptotically stable periodic solution that we
denote with $x_S(t)$, and therefore $\x_S=\bmat{x_S \\ x_S}$ will be a synchronous
solution of the coupled system, for all $\sigma\ge 0$.

\subsection{A piecewise smooth network}
The network equations are
\begin{equation}\begin{split}\label{Galvanetto_eq}
\dot y_1 = & \,\  y_2\\
\dot y_2 = & -y_1-\sigma(y_1-y_3)\pm\frac1{1+\gamma \abs{y_2-\bar v}}\\
\dot y_3 = & \,\ y_4\\
\dot y_4 = & -y_3-\sigma(y_3-y_1)\pm\frac1{1+\gamma \abs{y_4-\bar v}},
\end{split}
\end{equation}
with $\sigma, \gamma\ge 0$.  In the notation of \eqref{nonsmooth_eq}, we have
$$f^\pm(y)=\bmat{y_2 \\ -y_1\mp \frac{1}{1\pm \gamma (y_2-\bar v)}}, \,\ y \in \R^2\ , $$ 
and the discontinuity surface is the plane
$y_2 -\bar v=0$, 
so that $\nabla h=\bmat{0 \\1}$.

Further, with respect to the notation of \eqref{PWSnetwork}, here we have
$$x_1=\bmat{y_1\\  y_2},\quad x_2=\bmat{y_3\\ y_4}, \quad
\x=\bmat{x_1\\ x_2}, \,\ L=\left ( \begin{matrix} -1 & 1 \\ 1& -1 \end{matrix}\right),\,\  
E=\left ( \begin{matrix} 0 & 0 \\ 1& 0 \end{matrix}\right)$$ 
and rewrite system \eqref{Galvanetto_eq} as 
\begin{equation}\label{Galvanetto_2_eq}
	\dot \x=\left ( \begin{matrix} f^\pm(x_1) \\ f^\pm(x_2) \end{matrix} \right ) + 
	\sigma M\x, \quad M=L \otimes E .
\end{equation}

In the computations below, we fix $\bar v=0.15$, and
$\gamma=3$ as in \cite{Galvanetto_95}, and  use the MSF to
study the stability of the synchronous periodic solution as $\sigma$ varies in $[0,5]$. 
For $\sigma=0$, the two oscillators are uncoupled and the synchronous solution 
$\x_S=\bmat{x_S(t) \\ x_S(t)}$ 
has two Floquet multipliers at $1$ and two identical multipliers less than 1.
Of course, the synchronous solution 
persists for $\sigma \ge 0$, though its stability will depend on $\sigma$.
We use
the MSF to compute the Floquet exponents of the synchronous solution for $\sigma \neq 0$.  

Now we proceed like we did in Section \ref{MSF}.
Let $V=\begin{pmatrix} -1 & 1 \\1 & 1 \end{pmatrix}$, 
and $W=V \otimes I_2$, and let $S=W^{-1}MW$, so that
$S=\begin{pmatrix} -2E & {0} \\ {0} & {0} \end{pmatrix}$ where the $0$'s are 
${2 \times 2}$ blocks.  
Proceeding like we did to arrive at \eqref{Zi_eq}, we 
linearize \eqref{Galvanetto_2_eq} along $\x_S(t)$, 
and end up having to compute the Floquet multipliers of 
the linear non autonomous system
\begin{equation}\label{lin_single_eq}
\dot z=\left \{ \begin{matrix} (Df^\pm(x_S(t))-2\sigma E)  z, & {\nabla h^Tx_S(t)\gtrless \bar v},  \\
(Df_{\Sigma}(x_S(t))-2\sigma (E+B))  z, & {\nabla h^Tx_S(t)=\bar v},
\end{matrix} \right .
\end{equation} 
where $x_S(t)$ is the periodic solution of the single uncoupled agent $\dot x=f^\pm(x)$
and $B=\frac{(f^+-f^-)(x_S)\nabla h^T(x_S)}{\nabla h(x_S)^T(f^--f^+)(x_S)}E$.

Now, we make the observation that, for our problem \eqref{Galvanetto_2_eq},
we have $E+B=0$, and therefore in the sliding phase of \eqref{lin_single_eq} the
linearized problem is simply $\dot z=Df_{\Sigma}(x_S(t)) z$.

\subsection{Numerical experiments}
\begin{figure}\label{Mults}
		\includegraphics[width=250pt]{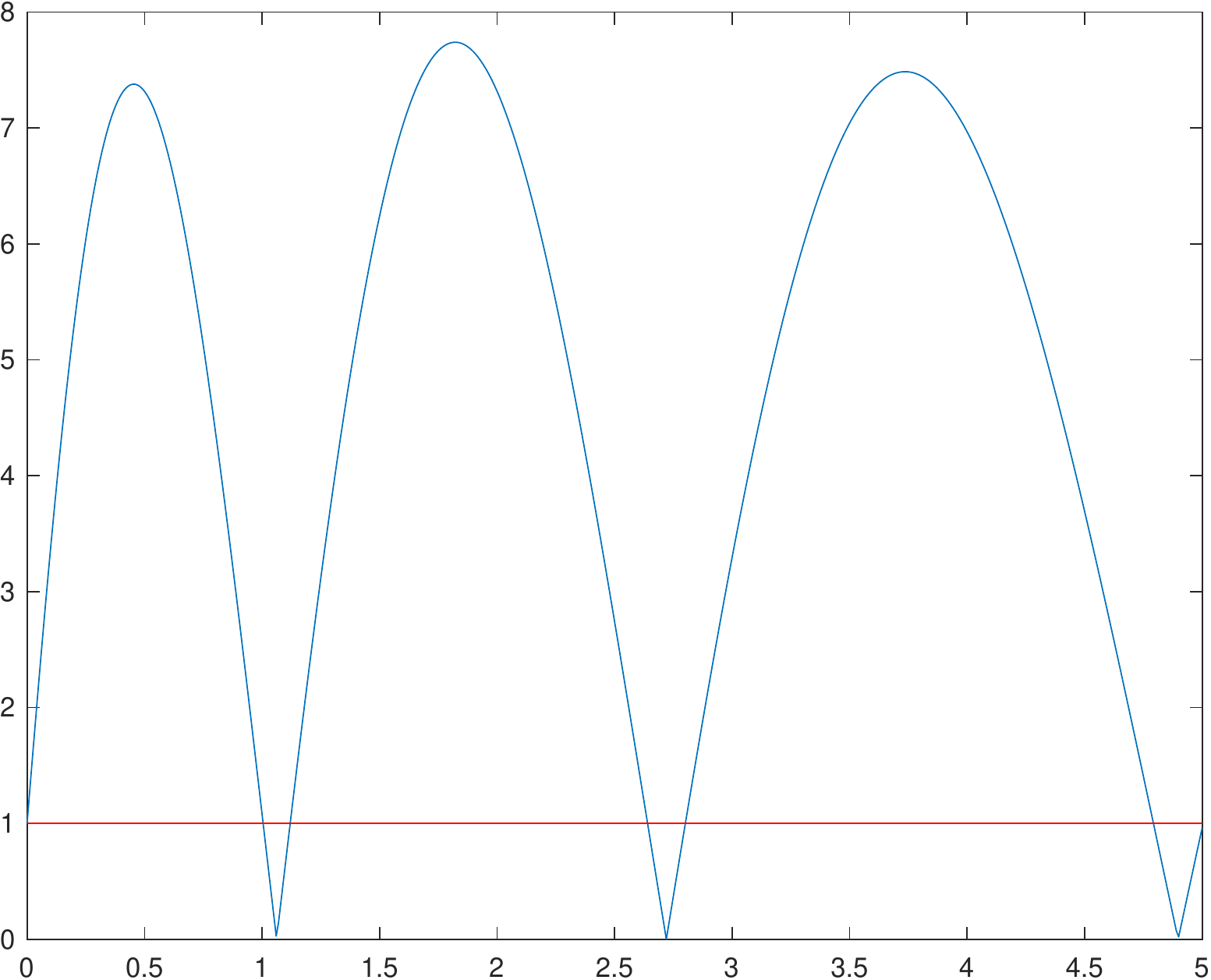}
\caption{Multiplier below 1 indicates stability of synchronized motion}
\end{figure}

Thanks to Theorem \ref{main_thm}, the Floquet exponents of the synchronous solution of 
\eqref{Galvanetto_eq} can be computed from \eqref{lin_single_eq}.   Therefore, our task
is to compute the solution of a single oscillator over one period, and then compute the
monodromy matrix of \eqref{lin_single_eq} and extract its Floquet multipliers, of which
we know one has to be $0$ because of sliding.  
Computation of the periodic orbit of the single agent is
done with the 4th order event technique of \cite{SlidNoProj} and fixed stepsize 
equal to $10^{-4}$ (so to
have a local error per step of size about $\mathtt{eps}$) and the monodromy matrix is
computed on the same mesh at once.  The multiplier at $0$ can always be 
recovered exactly due to the structure of the saltation matrix. 
In Figure \ref{Mults} we show the other multiplier in function of $\sigma$.
In the intervals where this is less than $1$, the synchronous solution is stable.

\begin{figure}\label{galvanetto_4p8_tol1em9_fig}
	$
	\begin{matrix}
		\includegraphics[width=200pt]{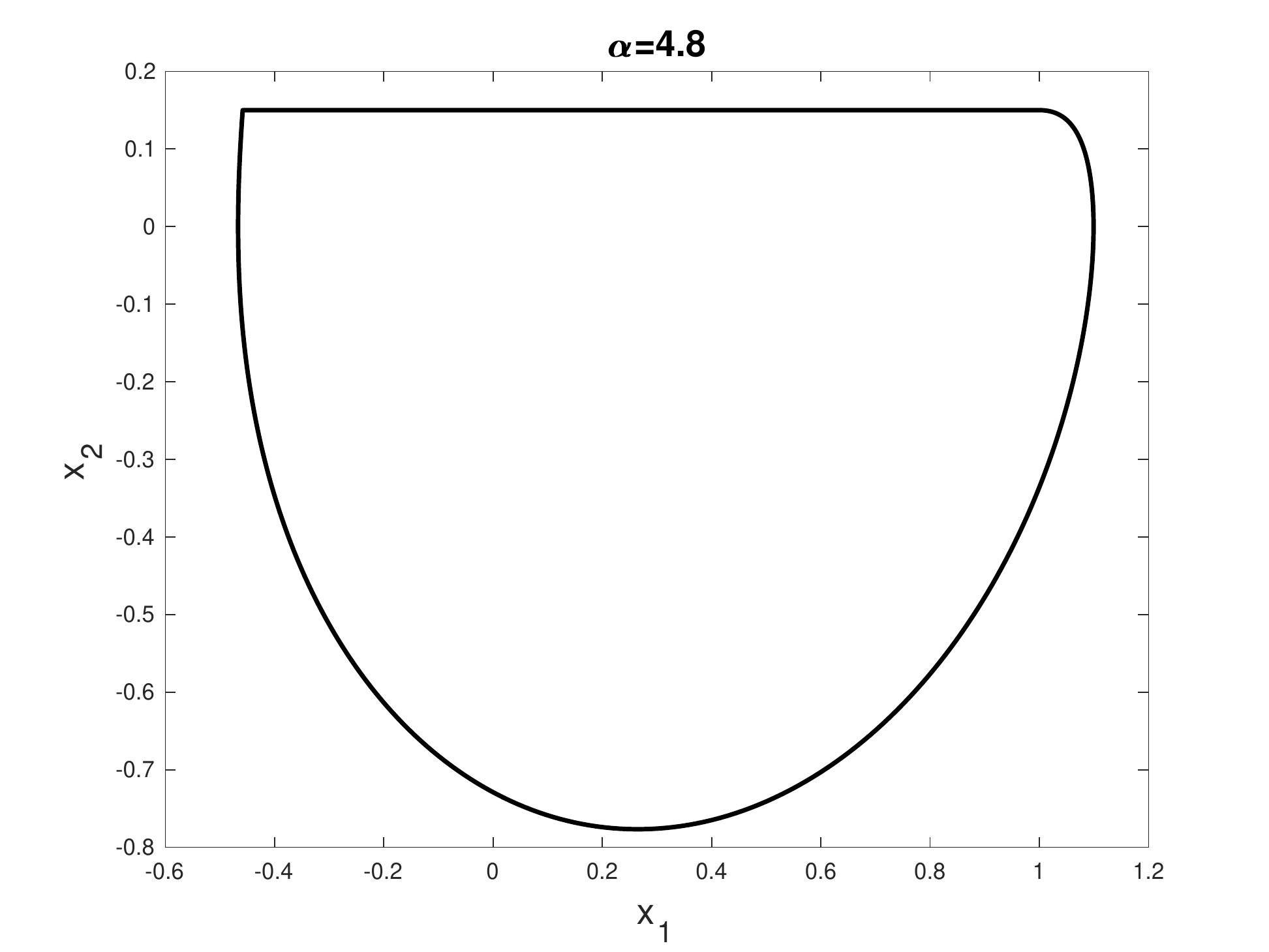} &
		\includegraphics[width=200pt]{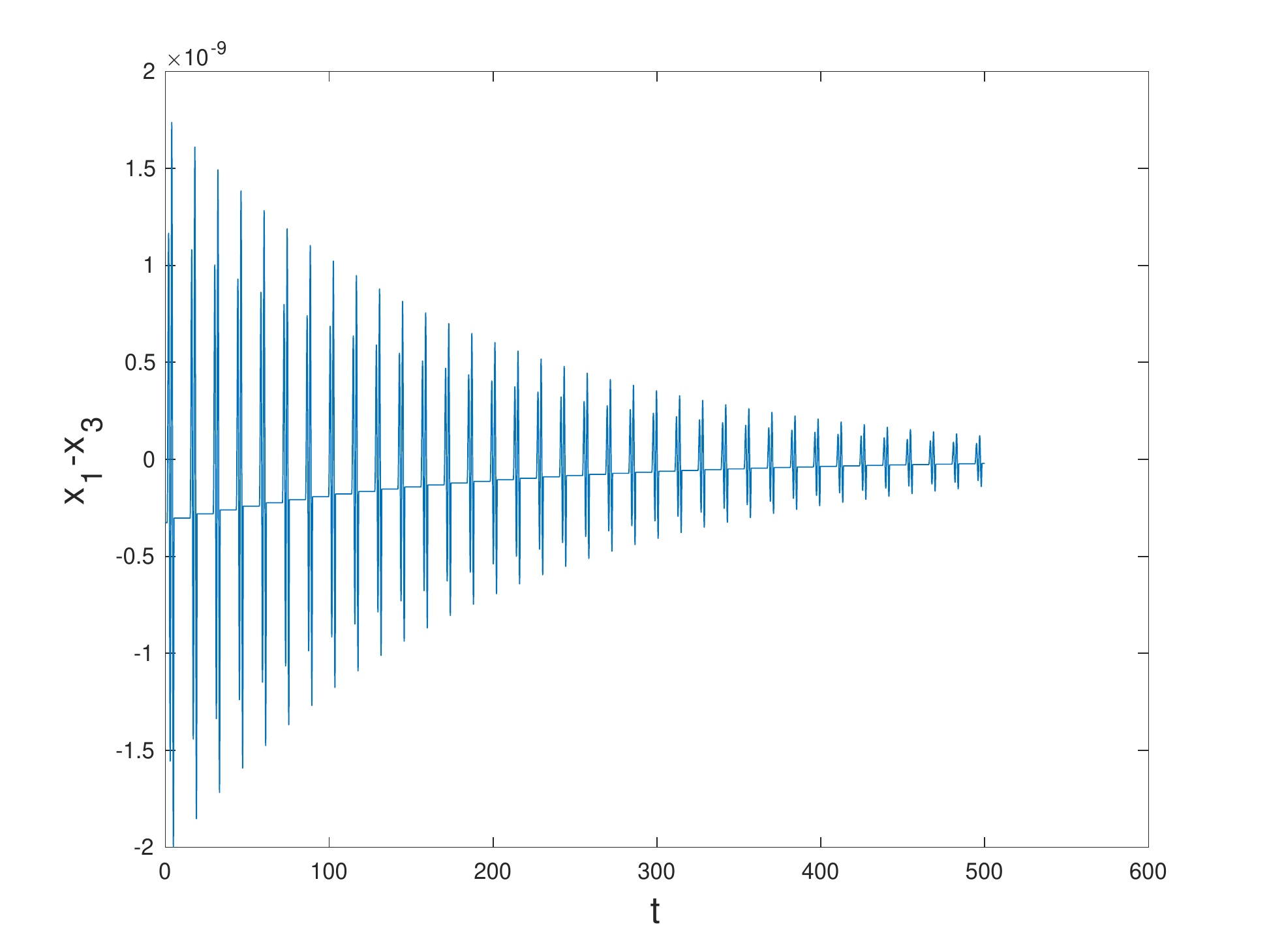}
	\end{matrix}
	$
\caption{Synchronous periodic orbit and network synchronization for $\sigma=4.8$}
\end{figure}

To confirm the results of the MSF analysis, 
we also integrated the full discontinuous system \eqref{Galvanetto_eq} with a variable 
stepsize integrator and event location techniques for sliding along the intersection of two discontinuity manifolds. For $\sigma=4.8$ the MSF plotted in Figure 
\ref{Mults} predicts asymptotic stability of the synchronous solution. 
Given initial conditions that do not belong to the synchronous 
manifold we integrated the full network for sufficiently large time to observe
convergence of the numerical solution to the synchronous periodic orbit. In Figure \ref{galvanetto_4p8_tol1em9_fig} on the left we plot the synchronous periodic
orbit while on the right we plot $(x_1(t)-x_3(t))$, after discarding the transient.  

We do not see synchronization for other parameter values such as 
$\sigma=1$, or $\sigma=2.6$, as predicted by the MSF, 
while we see synchronization for $\sigma=2.7$. 
Finally,
for the value of $\sigma=1.2$, our analysis based on the MSF validates the observation 
in \cite{Galvanetto_95} that the synchronous manifold is unstable.

\section{Conclusions}
In this work, we extended the Master Stability Function (MSF) tool to networks of
identical Piecewise Smooth  (PWS) differential systems,
in order to infer stability of a synchronous periodic solution of the
network.  Our analysis rested on the appropriate extension of the fundamental matrix
solution in the present PWS case.  We had to overcome several difficulties, in primis
the lack of uniqueness of suitable saltation matrices on the intersection of
several discontinuity manifolds and the possibility to decouple the (large) linearized
$nN$-system into $N$ systems of size $n$, in order to exploit the MSF technique.  We
succeeded in doing this under very general assumptions, for the network synchronizing
along a periodic orbit of a single agent.
We complemented our analysis by a numerical illustration of the use of the MSF
for a PWS system of mechanical oscillators synchronizing (for some values of
the coupling parameter) on a stick-slip oscillatory regime.
The case of synchronization on an orbit different from a periodic one remains
to be analyzed.

\bigskip
{\bf Acknowledgments} \\
This work has been partially supported by the GNCS-Indam group and the PRIN2017 research grant.  
The authors gratefully acknowledge the inspiration provided by a series of lectures given 
by Mario di Bernardo and Marco Coraggio at the University of Bari in March 2020. 
They were the last lectures before universities in Italy had to shut down due to the 
pandemic and provided a stimulating diversion in the following months.

\end{document}